\documentclass[11pt]{article}
\usepackage{amsthm,amsmath,latexsym,amssymb}

\newcommand{\bP}{{\rm |\kern-.15em P}}
\newcommand{\Q}{\kern.3em\rule{.07em}{.65em}\kern-.3em{\rm Q}}
\newcommand{\R}{{\rm I\kern-.15em R}}
\newcommand{\D}{{\rm |\kern-.15em D}}
\newcommand{\h}{{\rm |\kern-.15em H}}
\newcommand{\C}{\kern.3em\rule{.07em}{.65em}\kern-.3em{\rm C}}
\newcommand{\T}{{\rm T\kern-.35em T}}

\theoremstyle{plain}
\newtheorem{theorem}{Theorem}[section]
\newtheorem{lemma}[theorem]{Lemma}

\newtheorem{proposition}[theorem]{Proposition}
\newtheorem{corollary}[theorem]{Corollary}

\theoremstyle{definition}
\newtheorem{definition}[theorem]{Definition}
\newtheorem{example}[theorem]{Example}

\theoremstyle{remark}
\newtheorem{remark}[theorem]{Remark}

\begin{document}
\title{Injectivity of the composition operators of \'etale mappings}
\author{Ronen Peretz}
 
\maketitle

\begin{abstract}
Let $X$ be a topological space. The semigroup of all the \'etale mappings of $X$ (the local homeomorphisms $X\rightarrow X$), is denoted by ${\rm et}(X)$.
If $G\in {\rm et}(X)$ then the $G$-right (left) composition operator on ${\rm et}(X)$ is defined by:
$$
R_G\,\,(L_G):\,{\rm et}(X)\rightarrow {\rm et}(X),\,\,\,R_G(F)=F\circ G\,\,(L_G(F)=G\circ F).
$$
When are the composition operators injective? The first case we considered was 
that of entire functions $\mathbb{C}\rightarrow\mathbb{C}$ that are \'etale (and normalized). The second case
is that of the polynomial mappings $\mathbb{C}^2\rightarrow\mathbb{C}^2$ with determinant of their Jacobian matrix
equals (identically) to $1$ and whose $Y$-degrees equal their total degrees. For the first case we used the following definition:
$$
{\rm elh}(\mathbb{C})=\{ f\,:\,\mathbb{C}\rightarrow\mathbb{C}\,|\,f\,{\rm is}\,{\rm entire},\,
\forall\,z\in\mathbb{C}\,\,f'(z)\ne 0,\,\,f'(0)=1\}.
$$
Thus we use in this case the symbol ${\rm elh}(\mathbb{C})$ instead of ${\rm et}(X)$. Then we have: \\
\\
{\bf Proposition 3.12.} $\forall f\in {\rm elh}(\mathbb{C})$, $R_f$ in injective. \\
\\
{\bf Theorem 3.33. } Let $f\in {\rm elh}(\mathbb{C})$. Then $L_f$ is not injective if and only if
$$
f(z)=\frac{1}{b}e^{bz}+a\,\,{\rm for}\,\,{\rm some}\,\,a\in\mathbb{C},\,b\in\mathbb{C}^{\times}.
$$
These results that settled the first case were proved in: \\ 
R. Peretz, On the Structure of the Semigroup of Entire \'Etale Mappings,
{\em  The Journal Complex Analysis and Operator Theory}, Volume 7, Issue 5 (2013), Page 1655-1674. \\
The second case originated in a new approach to study \'etale polynomial
mappings $\mathbb{C}^2\rightarrow\mathbb{C}^2$ and in particular the two dimensional Jacobian Conjecture. This 
approach constructs a fractal structure on the semigroup of the (normalized) Keller mappings and outlines a new method of a possible attack
on this open problem (in preparation). The construction uses the left composition operator and the injectivity problem is essential. In this paper we will
completely solve the injectivity problems of the two composition operators for (normalized) Keller mappings. We will also solve the much
easier surjectivity problem of these composition operators.
\end{abstract}

\section{Introduction}
Let $X$ be a topological space. A mapping $F:\,X\rightarrow X$ is called a local homeomorphism of $X$, or an
\'etale mapping of $X$ if for any point $x\in X$ there exists a neighborhood $U$ of $x$ such that the restriction
of $F$ to $U$, denoted by $F|_U$, is an homeomorphism. The set of all the \'etale mappings of $X$, denoted by ${\rm et}(X)$,
is a semigroup with a unit with the composition of mappings taken to be the binary operation. If $G\in {\rm et}(X)$ then
the $G$-right composition operator on ${\rm et}(X)$ is defined by:
$$
\begin{array}{l} R_G:\,{\rm et}(X)\rightarrow {\rm et}(X) \\ R_G(F)=F\circ G.\end{array}
$$
The $G$-left composition operator on ${\rm et}(X)$ is defined by:
$$
\begin{array}{l} L_G:\,{\rm et}(X)\rightarrow {\rm et}(X) \\ L_G(F)=G\circ F.\end{array}
$$
We were interested in the injectivity of these two composition operators in two particular cases. The first is
the case of entire functions $\mathbb{C}\rightarrow\mathbb{C}$ that are \'etale (and normalized). The second case
is that of the polynomial mappings $\mathbb{C}^2\rightarrow\mathbb{C}^2$ with determinant of their Jacobian matrix
equals (identically) to $1$ and whose $Y$-degrees equal their total degrees. For the first case we use the following:
\begin{definition}
$$
{\rm elh}(\mathbb{C})=\{ f\,:\,\mathbb{C}\rightarrow\mathbb{C}\,|\,f\,{\rm is}\,{\rm entire},\,
\forall\,z\in\mathbb{C}\,\,f'(z)\ne 0,\,\,f'(0)=1\}.
$$
\end{definition}
\noindent
Thus we use in this case the symbol ${\rm elh}(\mathbb{C})$ instead of ${\rm et}(X)$. Then we have: \\
\\
{\bf Proposition 3.12. (in \cite{rp2})} $\forall f\in {\rm elh}(\mathbb{C})$, $R_f$ in injective. \\
\\
{\bf Theorem 3.33. (in \cite{rp2})} Let $f\in {\rm elh}(\mathbb{C})$. Then $L_f$ is not injective if and only if
$$
f(z)=\frac{1}{b}e^{bz}+a\,\,{\rm for}\,\,{\rm some}\,\,a\in\mathbb{C},\,b\in\mathbb{C}^{\times}.
$$
This settled the first case. It should be noted (see \cite{rp2}) that the proof for the left composition operator is much more involved
than the proof for the right composition operator (which follows directly from the Picard's Little Theorem). It is in fact the second case
that initiated our interest in the injectivity of the composition operators. It results from a new approach to study \'etale polynomial
mappings $\mathbb{C}^2\rightarrow\mathbb{C}^2$ and in particular the two dimensional Jacobian Conjecture (\cite{bcw},\cite{e} and \cite{s}). This 
approach constructs a fractal structure on the semigroup of the (normalized) Keller mappings and outlines a new method of a possible attack
on this open problem (in preparation). The construction uses the left composition operator and the injectivity problem is essential. In this paper we will
completely solve the injectivity problems of the two composition operators for (normalized) Keller mappings. We will also solve the much
easier surjectivity problem of these composition operators.

\section{The semigroup ${\rm et}(\mathbb{C}^2)$ of normalized Keller mappings and a few facts on their asymptotic variety}
Let $F\in\mathbb{C}[X,Y]^2$ be an \'etale mapping that satisfies the two normalizations: \\
1) $\det J_F\equiv 1$. \\
2) $\deg P=\deg_Y P$ and $\deg Q=\deg_Y Q$ where $F(X,Y)=(P(X,Y),Q(X,Y))\in\mathbb{C}[X,Y]^2$. \\
The set of all such mappings $F$ will be denoted by ${\rm et}(\mathbb{C}^2)$. This semigroup (with respect
to composition of mappings) is the parallel of the semigroup ${\rm elh}(\mathbb{C})$ for entire functions.
The $2$-dimensional Jacobian Conjecture can be rephrased in each of the following forms: \\
a) ${\rm et}(\mathbb{C}^2)\subseteq {\rm Aut}(\mathbb{C}^2)$. \\
b) $({\rm et}(\mathbb{C}^2),\circ)$ is a group. \\
For the next survey of results we refer to the following paper, \cite{rp5}.
We denote by $A(F)$ the asymptotic variety of $F$, i.e., the curve of all the asymptotic values of the mapping $F$. 
The canonical geometric basis of $F$ will be denoted by $R_0(F)$. This basis consists of finitely many rational 
mappings of the following form: $R(X,Y)=(X^{-\alpha},X^{\beta}Y+X^{-\alpha}\Phi(X))$, where $\alpha\in\mathbb{Z}^+$, 
$\beta\in\mathbb{Z}^+\cup\{ 0\}$, $\Phi(X)\in\mathbb{C}[X]$ and $\deg\Phi<\alpha+\beta$.
Also the effective $X$ powers in $X^{\alpha+\beta}Y+\Phi(X)$ have a gcd which equals $1$. The cardinality
of the geometric basis, $|R_0(F)|$, equals the number of components of the affine algebraic curve $A(F)$.
$\forall\,R\in R_0(F)$ we have the double asymptotic identity $F\circ R=G_R\in\mathbb{C}[X,Y]^2$ where the
polynomial mapping $G_R$ is called the $R$-dual of $F$. Each $R\in R_0(F)$ generates exactly one component of $A(F)$.
This component is normally parametrized by $\{G_R(0,Y)\,|\,Y\in\mathbb{C}\}$. We will denote by $H_R(X,Y)=0$
the implicit representation of this component in terms of the irreducible polynomial $H_R\in\mathbb{C}[X,Y]$. There
exists a natural number $\gamma(R)\ge 2$ and a polynomial $S_R(X,Y)\in\mathbb{C}[X,Y]$. The affine curve $S_R(X,Y)=0$ is 
called the $R$-phantom curve of $F$. The $R$-component of $A(F)$, $H_R(X,Y)=0$, is a polynomial curve which is not isomorphic 
to $\mathbb{A}^1$, and hence in particular must be a singular irreducible curve. We have the relation: $H_R(G_R(X,Y))=X^{\gamma(R)}S_R(X,Y)$.
The exponent $\gamma(R)$ satisfies the double inequality $2\le\gamma(R)\le\beta-\alpha$. In our case of the canonical 
rational mappings $R\in R_0(F)$, we have ${\rm sing}(R)=\{ X=0\}$. The following is true:
$$
G^{-1}_R(H_R(X,Y)=0)=G^{-1}_R(G_R({\rm sing}(R)))={\rm sing}(R)\cup\{S_R(X,Y)=0\}.
$$
Thus the $G_R$-preimage of the $R$-component of $A(F)$ (which is the $G_R$-image of ${\rm sing}(R)$) is
the union of two curves: the first is ${\rm sing}(R)$ and the second is the so-called $R$-phantom
curve of $F$. Even if for a single $R(X,Y)$ the $R$-phantom curve is empty then JC(2) follows. Also if $\forall\,R\in R_0(F)$
${\rm sing}(R)\cap \{ S_R(X,Y)=0\}=\emptyset$, then $F$ is a surjective mapping.

\begin{proposition}
If $F,G\in {\rm et}(\mathbb{C}^2)$ then $R_0(G)\subseteq R_0(F\circ G)$, $F(A(G))\subseteq A(F\circ G)$.
\end{proposition}
\noindent
{\bf Proof.} \\
$R\in R_0(G)\Rightarrow G\circ R\in\mathbb{C}[X,Y]^2\Rightarrow F\circ (G\circ R)\in\mathbb{C}[X,Y]^2
\Rightarrow (F\circ G)\circ R\in\mathbb{C}[X,Y]^2\Rightarrow R\in R_0(F\circ G)$. Next we have \\
$(a,b)\in F(A(G))\Rightarrow \exists\,R\in R_0(G)\exists\,Y\in \mathbb{C}\,\,{\rm such}\,{\rm that}\,\,
(a,b)=F((G\circ R)(0,Y))\Rightarrow\exists\,R\in R_0(F\circ G)\exists\,Y\in\mathbb{C}\,\,{\rm such}\,{\rm that}\,\,
(a,b)=((F\circ G)\circ R)(0,Y)\Rightarrow (a,b)\in A(F\circ G)$. $\qed $ \\
\\
The proposition tells us that compositions of \'etale mappings do not decrease the geometric basis of the 
right factor and consequently do not decrease the left image of its asymptotic variety. We naturally ask,
under what conditions the geometric basis of $F\circ G$ is actually larger than that of $G$? In other
words, we would like to know when is it true that $R_0(G)\subset R_0(F\circ G)$? This happens exactly
when $\exists\,R\in R_0(F\circ G)-R_0(G)$. This means that $(F\circ G)\circ R\in\mathbb{C}[X,Y]^2$,
$G\circ R\not\in\mathbb{C}[X,Y]^2$. Let $R(X,Y)=(X^{-\alpha},X^{\beta}Y+X^{-\alpha}\Phi(X))$,
$G(X,Y)=(P(X,Y),Q(X,Y))$. Then 
$$
(G\circ R)(X,Y)=(P(X^{-\alpha},X^{\beta}Y+X^{-\alpha}\Phi(X)),Q(X^{-\alpha},X^{\beta}Y+X^{-\alpha}\Phi(X))\in
$$
$$
\in\mathbb{C}(X,Y)^2-\mathbb{C}[X,Y]^2.
$$
We clearly have ${\rm sing}(G\circ R)\subseteq {\rm sing}(R)$ and so ${\rm sing}(G\circ R)=\{ X=0\}$. By 
$F\circ (G\circ R)=(F\circ G)\circ R\in\mathbb{C}[X,Y]^2$ we have $G\circ R\in R(F)$. This is not
necessarily a member of the geometric basis of $F$. The canonical geometric basis of $F$, $R_0(F)$ contains
finitely many rational mappings of the form: $S(X,Y)=(X^{-a},X^bY+X^{-a}\Psi(X))$.
Since $G\in {\rm et}(\mathbb{C}^2)$ it follows that $|\mathbb{C}^2-G(\mathbb{C}^2)|<\infty$ (a similar phenomenon
as the Picard's Little Theorem). If $L$ is an asymptotic tract of $F$ then $G^{-1}(L)$ can not be a bounded subset
of $\mathbb{C}^2$. The reason is that if $\overline{G^{-1}(L)}$ is compact, then $G(\overline{G^{-1}(L)})$ is
compact and since $L\subseteq G(G^{-1}(L))\subseteq G(\overline{G^{-1}(L)})$ this would imply the
contradiction that $L$ is bounded (and hence can not be an asymptotic tract). Hence $G^{-1}(L)$ has at least one
component, say $L_1$, that goes to infinity. This is because the number of components of $G^{-1}(L)$ is finite
and $G^{-1}(L)$ is not bounded. So $F\circ G$ has a limit along $L_1$ which equals the
above asymptotic value of $F$. This proves the following generalization of the second part of Proposition 2.1,
namely,

\begin{proposition}
If $F,G\in {\rm et}(\mathbb{C}^2)$ then $A(F)\cup F(A(G))=A(F\circ G)$.
\end{proposition}
\noindent
This proposition implies that if $A(F)\subset F(A(G))$ then necessarily $R_0(G)\subset R_0(F\circ G)$
because, as shown in the proof of Proposition 2.1 $\forall\,R\in R_0(G)$, $((F\circ G)\circ R)({\rm sing}(R))
\subseteq F(A(G))$.

\begin{proposition}
Let $F\in {\rm et}(\mathbb{C}^2)$. If $\exists\,G\in {\rm et}(\mathbb{C}^2)$ such that $R_0(G)=R_0(F\circ G)$,
then $F(\mathbb{C}^2)=\mathbb{C}^2$, i.e. $F$ is a surjective mapping.
\end{proposition}
\noindent
{\bf Proof.} \\
Since $F\in {\rm et}(\mathbb{C}^2)$ we have $\mathbb{C}^2-F(\mathbb{C}^2)\subseteq A(F)$, because in
this case the only points in the complement of the image of $F$ are the finitely many Picard exceptional 
values of $F$ which are asymptotic values of $F$. If, as the assumption says $R_0(G)=R_0(F\circ G)$
then by Proposition 2.2 we must have $A(F)\subseteq F(A(G))\subseteq F(\mathbb{C}^2)$ and so there are
no Picard exceptional values of the mapping $F$. $\qed $ \\
\\
If $F\in {\rm et}(\mathbb{C}^2)$ is not a surjective mapping then the last proposition implies that
$\forall\,G\in {\rm et}(\mathbb{C}^2)$ we must have $R_0(G)\subset R_0(F\circ G)$. In particular
$R_0(F)\subset R_0(F\circ F)$. This is the choice $G=F$. If we choose $G=F\circ F$ we get 
$R_0(F\circ F)\subset R_0(F\circ(F\circ F))$. Now it is clear that by induction we get the infinite
chain of strict inequalities:
$$
R_0(F)\subset R_0(F\circ F)\subset R_0(F\circ F\circ F)\subset\ldots\subset R_0(F^{\circ n})\subset\ldots
$$
where $F^{\circ n}=F\circ\ldots\circ F$ $n$ times fold composition. Since the cardinality of the geometric basis $R_0(G)$
is the number of components of the asymptotic variety $A(G)$ it follows in this case that the
asymptotic varieties of iterates of $F$ are of increasing complexity in the sense that the number of 
components of the curve $A(F^{\circ(n+1)})$ is strictly larger than the number of components of the 
curve $A(F^{\circ n})$. By Proposition 2.2 we always have $A(F)\cup F(A(F))= A(F\circ F)$. Hence
$A(F)\cup F(A(F))\cup (F\circ F)(A(F))=A(F)\cup F(A(F)\cup F(A(F)))=A(F)\cup F(A(F\circ F))=A(F\circ F\circ F)$.
By induction we get in general:
$$
A(F)\cup\bigcup^n_{k=1} F^{\circ K}(A(F))=A(F^{\circ(n+1)}).
$$
\section{The composition operators on ${\rm et}(\mathbb{C}^2)$ are not surjective but the right composition operator is injective}
\begin{proposition}
The mappings $R_F$, $L_F$ are not surjective if and only if $F\not\in {\rm Aut}(\mathbb{C}^2)$. In fact in this
case we have $R_F({\rm et}(\mathbb{C}^2))\subset {\rm et}(\mathbb{C}^2)-{\rm Aut}(\mathbb{C}^2)$,
$L_F({\rm et}(\mathbb{C}^2))\subset {\rm et}(\mathbb{C}^2)-{\rm Aut}(\mathbb{C}^2)$.
\end{proposition}
\noindent
{\bf Proof.} \\
$R_0(R_F(G))=R_0(G\circ F)\supseteq R_0(F)\ne\emptyset,\,\,\,A(L_F(G))=A(F\circ G)\supseteq A(F)\ne\emptyset$. $\qed $

\begin{proposition}
$R_F$ is injective.
\end{proposition}
\noindent
{\bf Proof.} \\
$R_F(G)=R_F(H)\Rightarrow G\circ F=H\circ F$. Since $F\in   {\rm et}(\mathbb{C}^2)$ we have $|\mathbb{C}^2-F(\mathbb{C}^2)|
<\infty$ and by the assumption $G|_{F(\mathbb{C}^2)}=H|_{F(\mathbb{C}^2)}$. Hence $G\equiv H$. $\qed $ \\
\\
We naturally inquire if also $L_F$ is injective. So let us assume that $L_F(G)=L_F(H)$. Then $F\circ G=F\circ H$. If we
denote $T=F\circ G$, then $R_0(G),\,R_0(H)\subseteq R_0(T)$ and also: $A(F)\cup F(A(G))=A(F)\cup F(A(H))$. If 
$A(F)=\emptyset$ then $F\in {\rm Aut}(\mathbb{C}^2)$ and so $G=F^{-1}(F\circ G)=F^{-1}(F\circ H)=H$. If, on the other hand,
$G\ne H$ then there are points $(X,Y)\in\mathbb{C}^2$ for which $G(X,Y)\ne H(X,Y)$. By our assumption:
$F(G(X,Y))=F(H(X,Y))$, so $F$ is not in ${\rm Aut}(\mathbb{C}^2)$ and it identifies different images of $G$ and $H$ of
the same $(X,Y)$. We ask the following question: Suppose that $F\circ G=F\circ H$, $G\ne H$. Is there a point
$(X,Y)\in\mathbb{C}^2$ for which $G(X,Y)=H(X,Y)$? \\
Based on our experience with entire functions we tend to prove that the answer to the question is
negative. Indeed this is the case and the proof is almost identical to the entire case, see \cite{rp2}. Namely, if the answer is 
affirmative, then we have two types of points in $\mathbb{C}^2$: those $(U,V)\in\mathbb{C}^2$ for which $G(U,V)\ne H(U,V)$
and the complimentary set, where both sets are non-empty. Let us denote by $N$ the first subset of $\mathbb{C}^2$, i.e.
$$
N=\{ (U,V)\in\mathbb{C}^2\,|\,G(U,V)\ne H(U,V)\}.
$$
The subset $N$ of $\mathbb{C}^2$ is open in the strong topology because $G$, $H$ are \'etale mappings and if 
$G(U,V)\ne H(U,V)$ then $\exists\,O$ an open neighborhood of $(U,V)$ in the strong topology such that 
$G(O)\cap H(O)=\emptyset$. So the complimentary subset of $N$ is a closed non-empty subset of $\mathbb{C}^2$.
Let $(X,Y)\in\partial N^c$ a boundary point of $N^c$. Let $(U_n,V_n)\in N$ satisfy $\lim (U_n,V_n)=(X,Y)$.
Then $\forall\,n\in\mathbb{Z}^+$, $G(U_n,V_n)\ne H(U_n,V_n)$, $F(G(U_n,V_n))=F(H(U_n,V_n))$, and $G(X,Y)=H(X,Y)$.
This implies that in any strong neighborhood of $G(X,Y)=H(X,Y)$ there are different points, say
$G(U_n,V_n)\ne H(U_n,V_n)$ for  $n\in\mathbb{Z}^+$ large enough, so that $F(G(U_n,V_n))=F(H(U_n,V_n))$. Hence
$F$ is not injective in any strong neighborhood of $G(X,Y)=H(X,Y)$. Thus $F\not\in {\rm et}(\mathbb{C}^2)$. This
contradiction proves the following,

\begin{proposition}
Suppose that $L_F$ is not injective for some $F\in {\rm et}(\mathbb{C}^2)$. Then $\forall\,G\ne H$,
$G,\,H\in {\rm et}(\mathbb{C}^2)$ such that $L_F(G)=L_F(H)$ and $\forall\,(X,Y)\in\mathbb{C}^2$ we have 
$G(X,Y)\ne H(X,Y)$.
\end{proposition}

\begin{remark}
Proposition 3.1 asserts that $\forall\,F\in {\rm et}(\mathbb{C}^2)-{\rm Aut}(\mathbb{C}^2)$ $\exists\,G,H
\in {\rm et}(\mathbb{C}^2)-{\rm Aut}(\mathbb{C}^2)$ such that $\forall\,M\in {\rm et}(\mathbb{C}^2)$ we
have $M\circ F\ne G$ and $F\circ M\ne H$.
\end{remark}
\section{The size, $d_F$, of the generic fiber of a Keller mapping $F\in{\rm et}(\mathbb{C}^2)$}
We will need the generic size of a fiber of a mapping $F=(P,Q)\in {\rm et}(\mathbb{C}^2)$. If we denote
$\deg P(X,Y)=n$ and $\deg Q(X,Y)=m$ then $\forall\,(a,b)\in\mathbb{C}^2$ the $F$ fiber over $(a,b)$ is
$F^{-1}(a,b)=\{(x,y)\in\mathbb{C}^2\,|\,F(x,y)=F(a,b)\}$. It is well known that this set is a finite subset
of $\mathbb{C}^2$ and, by the Bezout Theorem we have $|\{(x,y)\in\mathbb{C}^2\,|\,F(x,y)=(a,b)\}|=|F^{-1}(a,b)|\le n\cdot m$.
Moreover, there is a number that we will denote by $d_F$ such that generically in $(a,b)$ we have 
$|F^{-1}(a,b)|=d_F$. This means that $\{(a,b)\in\mathbb{C}^2\,|\,|F^{-1}(a,b)|\ne d_F\}$ is a closed
and proper Zariski subset of $\mathbb{C}^2$. In fact $\forall\,(a,b)\in\mathbb{C}^2,\,\,|F^{-1}(a,b)|\ne d_F
\Rightarrow |F^{-1}(a,b)|< d_F$. Thus we have $d_F=\max\{|F^{-1}(a,b)|\,|\,(a,b)\in\mathbb{C}^2\}$.

\begin{definition}
Let $F\in {\rm et}(\mathbb{C}^2)$. We will denote $d_F=\max\{|F^{-1}(a,b)|\,|\,(a,b)\in\mathbb{C}^2\}$. We will
call $d_F$ the geometrical degree of the \'etale mapping $F$.
\end{definition}

\begin{proposition}
$\forall\,F,G\in {\rm et}(\mathbb{C}^2),\,\,d_{F\circ G}=d_F\cdot d_G$.
\end{proposition}
\noindent
This is  a well known result. We include one of its proofs for convenience.
\\
{\bf Proof.} \\
$\forall\,(a,b)\in\mathbb{C}^2,\,\,(F\circ G)^{-1}(a,b)=G^{-1}(F^{-1}(a,b))$. But generically in $(a,b)$
$|F^{-1}(a,b)|=d_F$ and generically in $(c,d)$, $|G^{-1}(c,d)|=d_G$. $\qed $
\begin{definition}
An \'etale mapping $F\in {\rm et}(\mathbb{C}^2)$ is composite if $\exists\,G,H\in {\rm et}(\mathbb{C}^2)-
{\rm Aut}(\mathbb{C}^2)$ such that $F=G\circ H$. An \'etale mapping $A\in {\rm et}(\mathbb{C}^2)-{\rm Aut}(\mathbb{C}^2)$ is
prime if it is not composite. This is equivalent to: $A=B\circ C$ for some $B,C\in  {\rm et}(\mathbb{C}^2)$
$\Rightarrow B\in {\rm Aut}(\mathbb{C}^2)\vee C\in {\rm Aut}(\mathbb{C}^2)$. The subset of ${\rm et}(\mathbb{C}^2)$
of all the prime mappings will be denoted by ${\rm et_p}(\mathbb{C}^2)$. Thus the set of all the composite
\'etale mappings is ${\rm et}(\mathbb{C}^2)-{\rm et_p}(\mathbb{C}^2)$.
\end{definition}

\begin{proposition}
$\forall\,F\in {\rm et}(\mathbb{C}^2)-{\rm et_p}(\mathbb{C}^2)$, $d_F$ is not a prime number. Equivalently,
$\forall\,F\in {\rm et}(\mathbb{C}^2)$, $d_F$ is a prime number $\Rightarrow F\in {\rm et_p}(\mathbb{C}^2)$.
\end{proposition}
\noindent
{\bf Proof.} \\
$F\in {\rm et}(\mathbb{C}^2)-{\rm et_p}(\mathbb{C}^2)\Rightarrow\exists\,G,H\in {\rm et}(\mathbb{C}^2)-
{\rm Aut}(\mathbb{C}^2)$ such that $F=G\circ H$ (by the definition) $\Rightarrow d_F=d_G\cdot d_H,\,\,d_G,d_H>1$
(by Proposition 4.2 and the fact $d_M=1\Leftrightarrow M\in {\rm Aut}(\mathbb{C}^2)$) $\Rightarrow d_F$ is
a composite integer. $\qed $

\begin{theorem} The following hold true: \\
{\rm 1)} If $ {\rm et}(\mathbb{C}^2)-{\rm Aut}(\mathbb{C}^2)\ne\emptyset$ then ${\rm et_p}(\mathbb{C}^2)\ne\emptyset$ \\
{\rm 2)} $\forall\,F\in {\rm et}(\mathbb{C}^2)\,\,\exists\,k\in\mathbb{Z}^+\cup\{ 0\}\,\,\exists\,
A_0\in {\rm Aut}(\mathbb{C}^2)\,\exists\,P_1,\ldots,P_k\in {\rm et_p}(\mathbb{C}^2)$ such that
$F=A_0\circ P_1\circ\ldots\circ P_k$.
\end{theorem}
\noindent
{\bf Proof.} \\
If ${\rm et_p}(\mathbb{C}^2)=\emptyset $ then ${\rm et}(\mathbb{C}^2)-{\rm Aut}(\mathbb{C}^2)$ are all
composite \'etale mappings. Let $F\in {\rm et}(\mathbb{C}^2)-{\rm Aut}(\mathbb{C}^2)$, then
$\exists\,G_1,G'_2\in {\rm et}(\mathbb{C}^2)-{\rm Aut}(\mathbb{C}^2)$ such that $F=G_1\circ G'_2$. So
$\exists\,G_2,G'_3\in {\rm et}(\mathbb{C}^2)-{\rm Aut}(\mathbb{C}^2)$ such that $G'_2=G_2\circ G'_3$.
Hence $F=G_1\circ G_2\circ G'_3$. Continuing this we get for any $k\in\mathbb{Z}^+$ $\exists\,G_1,\ldots,G_k\in
{\rm et}(\mathbb{C}^2)-{\rm Aut}(\mathbb{C}^2)$ such that $F=G_1\circ\ldots\circ G_k$ and by Proposition 4.2
$d_F=\prod_{j=1}^k d_{G_j}$. But $\forall\,1\le j\le k$, $d_{G_j}\ge 2$ and so $\forall\,k\in\mathbb{Z}^+$,
$d_F\ge 2^k$ a contradiction to $d_F<\infty $. Thus ${\rm et_p}(\mathbb{C}^2)\ne\emptyset$. \\
Now part 2 is standard, for if $F\in {\rm Aut}(\mathbb{C}^2)$ we take $A_0=F$ and $k=0$. If $F\in
{\rm et_p}(\mathbb{C}^2)$ we take $A_0={\rm id}$, $k=1$, and $P_1=F$. If $F\in {\rm et}(\mathbb{C}^2)-
{\rm et_p}(\mathbb{C}^2)$ then $F=G\circ H$ for some $G,H\in {\rm et}(\mathbb{C}^2)-{\rm Aut}(\mathbb{C}^2)$.
So by Proposition 4.2 $d_F=d_G\cdot d_H$ and since $d_G,d_H\ge 2$ it follows that $d_G,d_H<d_F$ and
we conclude the proof of part 2 using induction on the geometrical degree. Namely $G=P_1\circ\ldots\circ P_m$,
$H=P_{m+1}\circ\ldots\circ P_k$ for $m\ge 1$, $k\ge m+1$ and some primes $P_1,\ldots,P_k\in {\rm et_p}(\mathbb{C}^2)$ $\qed $

\section{The metric spaces $({\rm et}(\mathbb{C}^2),\rho_D)$}
We will need a special kind of four (real) dimensional subsets of $\mathbb{R}^4$. These will serve us to construct suitable metric structures 
on ${\rm et}(\mathbb{C}^2)$.

\begin{definition}
Let $D$ be an open subset of $\mathbb{C}^2$ with respect to the strong topology, that satisfies the following conditions: \\
1) ${\rm int}(\overline{D})=D$ ($D$ has no "slits"). \\
2) $\overline{D}$ is a compact subset of $\mathbb{C}^2$ (in the strong topology). \\
3) $\forall\,G_1,G_2\in {\rm et}(\mathbb{C}^2)$, $G_1(D)=G_2(D)\Leftrightarrow G_1=G_2$. \\
We define the following real valued function:
$$
\rho_D\,:\,{\rm et}(\mathbb{C}^2)\times {\rm et}(\mathbb{C}^2)\rightarrow\mathbb{R}^+\cup\{0\},
$$
$$
\rho_D(G_1,G_2)={\rm the\,volume\,of}\,\,G_1(D)\Delta G_2(D).
$$
Here we use the standard set-theoretic notation of the symmetric difference between two sets $A$ and $B$, i.e.
$A\Delta B=(A-B)\cup (B-A)$.
\end{definition}

\begin{remark}
It is not clear how to construct an open subset $D$ of $\mathbb{C}^2$ that will satisfy the three properties
that are required in definition 5.1. We will postpone for a while the demonstration that such open sets exist.
\end{remark}

\begin{proposition}
$\rho_D$ is a metric on ${\rm et}(\mathbb{C}^2)$.
\end{proposition}
\noindent
{\bf Proof.} \\
1) $\rho_D(G_1,G_2)=0 \Leftrightarrow {\rm the\,volume\,of}\,\,G_1(D)\Delta G_2(D)=0 \Leftrightarrow G_1(D)=G_2(D)$
(where the last equivalence follows by the fact that $G_1$ and $G_2$ are local homeomorphisms in the strong topology and because
of condition 1 in definition 5.1) $\Leftrightarrow G_1=G_2$ (by condition 3 in definition 5.1). \\
2) By $G_1(D)\Delta G_2(D)=G_2(D)\Delta G_1(D)$ it follows that $\rho_D(G_1,G_2)=\rho_D(G_2,G_1)$. \\
3) Here we use a little technical set-theoretic containment. Namely, for any three sets $A,B$ and $C$ we have, 
$$
A\Delta C\subseteq (A\Delta B)\cup (B\Delta C).
$$
This implies that $G_1(D)\Delta G_3(D)\subseteq (G_1(D)\Delta G_2(D))\cup (G_2(D)\Delta G_3(D))$ from which it follows that
$$
({\rm the\,volume\,of}\,\,G_1(D)\Delta G_3(D))\le ({\rm the\,volume\,of}\,\,G_1(D)\Delta G_2(D))+
$$
$$
+({\rm the\,volume\,of}\,\,G_2(D)\Delta G_3(D)).
$$
Hence the triangle inequality $\rho_D(G_1,G_3)\le\rho_D(G_1,G_2)+\rho_D(G_2,G_3)$ holds. $\qed $ \\
\\
So far we thought of the volume of $G_1(D)\Delta G_2(D)$ as the volume of the open set which is the symmetric difference between
the $G_1$ image and the $G_2$ image of the open set $D$. However, the mappings $G_1$ and $G_2$ are \'etale and in particular need
not be injective. We will take into the volume computation the multiplicities of $G_1$ and of $G_2$. By Theorem 3 on page 39 of \cite{bm}
we have the following: Given $F\,:\,\mathbb{C}^n\rightarrow\mathbb{C}^n$ we define $\tilde{F}=({\rm Re}\,F_1,{\rm Im}\,F_1,\ldots,
{\rm Re}\,F_n,{\rm Im}\,F_n)\,;\,\mathbb{R}^{2n}\rightarrow\mathbb{R}^{2n}$. Then $\det J_{\tilde{F}}=|\det J_F|^2$. Thus the Jacobian 
Condition, $\det J_F\equiv 1$ implies that $\det J_{\tilde{F}}\equiv 1$. So the real mapping $\tilde{F}$ preserves the usual volume form.
In order to take into account the multiplicities of the \'etale mappings $G_1$ and $G_2$ when computing the volume of the symmetric
difference $G_1(D)\Delta G_2(D)$ we had to do the following. For any $G\in {\rm et}(\mathbb{C}^2)$ instead of computing,
$$
\int\int\int\int_D\left(\det J_{\tilde{G}}\cdot dV\right)=\int\int\int\int_{D}dV,
$$
we compute
$$
\int\int\int\int_{\tilde{G}(D)}dX_1 dX_2 dY_1 dY_2\,\,\,\,\,{\rm where}\,\,X=X_1+iX_2,\,Y=Y_1+iY_2.
$$
For every $j=1,2,\ldots,d_G$ we denote by $D_j$ that subset of $D$ such that for each point of $D_j$ there are exactly $j$ points 
of $D$ that are mapped by $G$ to the same image of that point. In other words, $D_j=\{\alpha\in D\,|\,|\tilde{G}^{-1}(\tilde{G}(\alpha))\cap D|=j\}$.
We assume that $D$ is large enough so that $\forall\,j=1,\ldots,d_G$ we have $D_j\ne\emptyset$.
For our \'etale mappings it is well known that if $j<d_G$ then $\dim D_j<\dim D$ so the volume these $D_j$'s contribution equals to 0. 
The dimension claim follows by the well known fact that the size of a generic fiber $|G^{-1}(x)|$ equals to $d_G$ and that $d_G$ is also the
maximal size of any of the fibers of $G$.
However, for the sake of treating more general families of mappings we denote by ${\rm vol}(D_j)$ the volume of the set $D_j$.
Then $D$ has a partition into exactly $j$ subsets of equal volume. The volume of each such a set is ${\rm vol}(D_j)/j$ and each such a set has
exactly one of the $j$ points in $\tilde{G}^{-1}(\tilde{G}(\alpha))\cap D$ for each $\alpha\in D_j$. We note that ${\rm vol}(\tilde{G}(D_j))={\rm vol}(D_j)/j$
by the Jacobian Condition. Thus the volume with the multiplicity of $\tilde{G}$ taken into account is given by:
$$
{\rm vol}(\tilde{G}(D))+\sum_{j=2}^{d_G}(j-1)\cdot\frac{{\rm vol}(D_j)}{j}=
{\rm vol}(\tilde{G}(D))+\sum_{j=2}^{d_G}(j-1)\cdot{\rm vol}(\tilde{G}(D_j)).
$$
We note that $\tilde{G}(D)=\bigcup_{j=1}^{d_G}\tilde{G}(D_j)$ is a partition, so ${\rm vol}(\tilde{G}(D))=\sum_{j=1}^{d_G}{\rm vol}(\tilde{G}(D_j))$. 
Hence we can express the desired volume by
$$
{\rm vol}(\tilde{G}(D))+\sum_{j=2}^{d_G}(j-1)\cdot{\rm vol}(\tilde{G}(D_j))=\sum_{j=1}^{d_G}j\cdot {\rm vol}(\tilde{G}(D_j)).
$$
We note that this equals to $\sum_{j=1}^{d_G}{\rm vol}(D_j)$ and since $D=\bigcup_{j=1}^{d_G} D_j$ is a partition we have
${\rm vol}(D)=\sum_{j=1}^{d_G} {\rm vol}(D_j)$. As expected, the volume computation that takes into account the multiplicity of $G$ is in general
larger than the geometric volume ${\rm vol}(\tilde{G}(D))$. The access can be expressed in several forms:
$$
{\rm vol}(D)-{\rm vol}(\tilde{G}(D))=\sum_{j=2}^{d_G}(j-1)\cdot {\rm vol}(\tilde{G}(D_j))=\sum_{j=2}^{d_G}\left(1-\frac{1}{j}\right){\rm vol}(D_j).
$$
Coming back to the computation of the metric distance $\rho_D(G_1,G_2)={\rm the\,volume\,of}\,\,G_1(D)\Delta G_2(D)$ we compute the volume
of $G_1(D)-G_2(D)$ with the multiplicity of $G_1$ while the volume of $G_2(D)-G_1(D)$ is computed with the multiplicity of $G_2$.

\section{Characteristic sets of families of holomorphic local homeomorphisms $\mathbb{C}^2\rightarrow\mathbb{C}^2$}
In this section we prove the existence of sets $D$ that satisfy the three properties that are required in definition 5.1.
The third property will turn out to be the tricky one. 

\begin{definition}
let $\Gamma$ be a family of holomorphic local homeomorphisms $F\,:\,\mathbb{C}^2\rightarrow\mathbb{C}^2$.
A subset $D\subseteq\mathbb{C}^2$ is called a characteristic set of $\Gamma$ if it satisfies the
following condition: $\forall\,F_1,F_2\in\Gamma$, $F_1(D)=F_2(D)\Leftrightarrow F_1=F_2$.
\end{definition}
\noindent
We start by recalling the well known rigidity property of holomorphic functions in one complex variable. Also known as the permanence principle,
or the identity theorem. The identity theorem  for analytic functions of one complex variable says that if $D\subseteq\mathbb{C}$ is a domain
(an open and a connected set) and if $E$ is a subset of $D$ that has a non-isolated point and if $f(z)$ is an analytic function defined on $D$ and 
vanishing on $E$, then $f(z)=0$ for all $z\in D$.

There is an identity theorem for analytic functions in several complex variables, but for more than one variable the above statement is false.
One correct statement is as follows:"Let $f(z)$ and $g(z)$ be holomorphic functions in a domain $D$ in $\mathbb{C}^n$. If $f(z)=g(z)$ for all $z$ 
in a non-empty set $\delta$ in $D$, then $f(z)=g(z)$ in $D$. Hence, analytic continuation of holomorphic functions in several complex variables
can be performed as in the case of one complex variable. Contrary to the case of one complex variable, the zero set of a holomorphic function in 
a domain $D\subseteq\mathbb{C}^n$, $n\ge 2$, contains no isolated points. Thus even if $f(z)=g(z)$ in a set with accumulation points in $D$, it
does not necessarily follow that $f(z)=g(z)$ in $D$. For example, in $\mathbb{C}^2$ with variables $z$ and $w$ we can take $f(z,w)=z$ and
$g(z,w)=z^2$." (Chapter 1, page 16 in \cite{tn}).

In spite of the above standard identity theorem for $n\ge 2$ complex variables, that requires a thick set $E$ (i.e. an open set) on which $f(z)=g(z)$
one can do much better. Let us start with the following. Let $F(Z,W)$ be an entire function of two complex variables $Z$ and $W$. Let us define a
subset $E$ of $\mathbb{C}^2$ as follows. We take a convergent sequence $\{Z_k\}_{k=1}^{\infty}$ of different numbers. Thus $\lim Z_k=a$ and
$j\ne k\Rightarrow Z_j\ne Z_k$. For each $k$, let $\{W_j^{(k)}\}_{j=1}^{\infty}$ be a convergent sequence of different numbers, such that
their limit is $\lim_{j\rightarrow\infty} W_j^{(k)}=Z_k^{'}$. We define $E=\{(Z_k,W_j^{(k)})\,|\,j,k=1,2,3,\ldots\}$. Now we have,

\begin{proposition}
If $F(Z,W)$ vanishes on $E$, i.e. $F(Z_k,W_j^{(k)})=0$ for $j,k=1,2,3,\ldots $, then $F(Z,W)\equiv 0$ is the zero function.
\end{proposition}

\begin{remark}
We note that $E$ is a thin set, in fact a countable set. Even the closure $\overline{E}$ is thin.
\end{remark}
\noindent
{\bf A Proof of Proposition 6.2.} \\
Since $F(Z,W)$ is an entire function, it can be represented as a convergent power series centered at $(0,0)$ with an infinite
radius of convergence. We can sum the terms in the order we please. Let us write $F(Z,W)$ as a power series in $W$ with
coefficients that are entire functions in $Z$. Thus we have, $F(Z,W)=\sum_{k=0}^{\infty} a_k(Z)W^k$, where for each $k=0,1,2,\ldots $,
$a_k(Z)$ is an entire function in the variable $Z$. For a fixed $l\in\{1,2,3,\ldots\}$ we have by our assumptions the following,
$F(Z_l,W_j^{(l)})=0$ for $j=1,2,3,\ldots $. But $\lim_{j\rightarrow\infty}W_j^{(l)}=Z_l^{'}$ so that $g_l(W)=F(Z_l,W)$ is an
entire function of the single variable $W$, which vanishes on a convergent sequence $\{W_j^{(l)}\}_{j=1}^{\infty}$. By the identity
theorem of one complex variable we deduce that $g_l(W)\equiv 0$, the zero function. Since, $g_l(W)=\sum_{k=0}^{\infty} a_k(Z_l)W^k$ it
follows that the Maclaurin coefficients $a_k(Z_l)$, $k=0,1,2,\ldots $ vanish. Now, this is valid for each $l$, and $\lim Z_l=a$ converges.
Since each $a_k(Z)$ is an entire function which vanishes on a convergent sequence $\{ Z_l\}_{l=1}^{\infty}$ it follows, once again, by
the identity theorem in one complex variable, that $a_k(Z)\equiv 0$, the zero function, $k=0,1,2,\ldots $. Hence we conclude that
$F(Z,W)=\sum_{k=0}^{\infty} a_k(Z)W^k\equiv 0$. $\qed $ \\
\\
This type of elementary arguments that was used to construct a thin set $E$ for identity purpose, is not new. For example: \\
"Theorem. Let $D\subseteq\mathbb{C}$ be a domain, and let $E$ be a subset of $D$ that has a non-isolated point. Let $F(Z,W)$ be a function 
defined for $Z,W\in D$ such that $F(Z,W)$ is analytic in $Z$ for each fixed $W\in D$ and analytic in $W$ for each fixed $Z\in D$.
If $F(Z,W)=0$ whenever $Z$ and $W$ both belong to $E$, then $F(Z,W)=0$ for all $Z,W\in D$.", \cite{g}.

Advancing along the lines of the construction of the thin set in Proposition 6.2 we note that if $\{Z_k\}_{k=1}^{\infty}$ is a
sequence of different numbers that converges to $\lim Z_k=a$, and if for each $k=1,2,3,\ldots $ there is a straight line segment 
$[\alpha_k,\beta_k]$ of $W$'s such that two entire functions $F(Z,W)$ and $G(Z,W)$ agree on the union (a countable union) of
the segments $\{Z_k\}\times [\alpha_k,\beta_k]$, i.e. $F(Z_k,W)=G(Z_k,W)$, $\forall\,W\in[\alpha_k,\beta_k]$, then $F(Z,W)\equiv G(Z,W)$,
$\forall\,(Z,W)\in\mathbb{C}^2$.

We now will construct characteristic sets of families $\Gamma$ of holomorphic local homeomorphisms $F:\,\mathbb{C}^2\rightarrow\mathbb{C}^2$.

\begin{definition}
Let $m$ be a natural number and $\alpha\in\mathbb{C}^2$. An $m$-star at $\alpha$ is the union of $m$ line segments, so that any pair intersect
in $\alpha $.
\end{definition}

\begin{definition}
Let $l$ be a line segment and let $\{\alpha_k\}$ be a countable dense subset of $l$. Let $\{n_k\}$ be a sequence of different natural numbers
and $\forall\,k$, let $S_{n_k}$ be an $n_k$-star at $\alpha_k$ such that one of the star's segments lies on $l$, and such that $\forall\,k_1\ne k_2$,
$\tilde{S}_{n_{k_1}}\cap\tilde{S}_{n_{k_2}}=\emptyset$. Here we denoted $\tilde{S}=S-l$. Moreover, we group the stars in bundles of, say 5, thus
getting the sequence of star bundles:
$$
\{S_{n_1},S_{n_2},S_{n_3},S_{n_4},S_{n_5}\},\,\{S_{n_6},\ldots,S_{n_{10}}\},\ldots,\{S_{n_{5j+1}},\ldots,S_{n_{5j+5}}\},\ldots
$$
and for each bundle of five we take the maximal length of its rays to be at most $1/10$ the length of the maximal length of the previous
bundle. We define,
$$
l_0^{\{n_k\}}=l\cup\bigcup_{k=1}^{\infty}S_{n_k}.
$$
Let $\{Z_k\}_{k=1}^{\infty}$ be a sequence of different complex numbers that converges to $\lim Z_k=a$. Let $\{\{n_j^{(k)}\}_{j=1}^{\infty}\}_{k=1}^{\infty}$
be a partition of the natural numbers, $\mathbb{Z}^+$. In fact all we need is the disjointness, i.e. $k_1\ne k_2\Rightarrow
\{n_j^{(k_1)}\}_{j=1}^{\infty}\cap\{n_j^{(k_2)}\}_{j=1}^{\infty}=\emptyset$. Let us consider the stared segments
$$
\{l_0^{\{n_j^{(k)}\}_{j=1}^{\infty}}\,|\,k=1,2,3,\ldots\}
$$
and define the following countable union of stared segments in $\mathbb{C}^2$:
$$
\bigcup_{k=1}^{\infty}\{Z_k\}\times l_0^{\{n_j^{(k)}\}_{j=1}^{\infty}},
$$
where we assume that the lengths of the star rays were chosen to satisfy disjointness in $\mathbb{C}^2$, namely:
$$
k_1\ne k_2\Rightarrow\{Z_{k_1}\}\times l_0^{\{n_j^{(k_1)}\}_{j=1}^{\infty}}\cap\{Z_{k_2}\}\times l_0^{\{n_j^{(k_2)}\}_{j=1}^{\infty}}=\emptyset.
$$
We let,
$$
E=\bigcup_{k=1}^{\infty}\{Z_k\}\times l_0^{\{n_j^{(k)}\}_{j=1}^{\infty}},
$$
or if we need a closed (compact) set, the closure of this union.
\end{definition}

\begin{proposition}
Let $\Gamma$ be any family of entire holomorphic local homeomorphisms $F:\,\mathbb{C}^2\rightarrow\mathbb{C}^2$. Then $E$ is a characteristic 
set of $\Gamma$.
\end{proposition}
\noindent
{\bf Proof.} \\
Let $F_1,F_2\in\Gamma$ satisfy $F_1(E)=F_2(E)$. Then each stared line segment,
$$
\{Z_k\}\times l_0^{\{n_j^{(k)}\}_{j=1}^{\infty}},
$$
must be mapped onto a curve,
$$
F_1(\{Z_k\}\times l_0^{\{n_j^{(k)}\}_{j=1}^{\infty}})=F_2(\{Z_k\}\times l_0^{\{n_j^{(k)}\}_{j=1}^{\infty}})
$$
and each $n_j^{(k)}$-star on $l$, $S_{n_j^{(k)}}$ is mapped onto a holomorphic $n_j^{(k)}$-star,
$$
F_1(\{Z_k\}\times S_{n_j^{(k)}})=F_2(\{Z_k\}\times S_{n_j^{(k)}}).
$$
This is because the valence sequences of the stars 
$$
\{\{n_j^{(k)}\}_{j=1}^{\infty}\}_{k=1}^{\infty},
$$
are pairwise disjoint natural numbers, and $F_1$, $F_2$ are local homeomorphisms and hence preserve the star
valencies $n_j^{(k)}$. The centers of the holomorphic stars,
$$
\{F_1(\alpha_{n_j^{(k)}})\}=\{F_2(\alpha_{n_j^{(k)}})\},
$$
form a countable and a dense subset of the curves $F_1(\{Z_k\}\times l)=F_2(\{Z_k\}\times l)$. By continuity this
implies that the restrictions,
$$
F_1|_{\{Z_k\}\times l}\,\,\,\,\,{\rm and}\,\,\,\,\,F_2|_{\{Z_k\}\times l},
$$
coincide. Since $F_1$ and $F_2$ are holomorphic, this implies by Proposition 6.2 (which is a variant of the identity
theorem for entire functions $\mathbb{C}^2\rightarrow\mathbb{C}^2$) that $F_1\equiv F_2$. $\qed $ \\

\begin{remark}
Proposition 6.6 holds true for any rigid family of local homeomorphisms. Rigidity here means that 
$$
F_1|_{\{Z_k\}\times l}=F_2|_{\{Z_k\}\times l}\Leftrightarrow F_1\equiv F_2.
$$ 
So the proposition holds true for holomorphic mappings, for harmonic mappings and in particular for ${\rm et}(\mathbb{C}^2)$.
\end{remark}
\noindent
We recall that definition 5.1 required also two additional topological properties, namely the open set $D$ should satisfy ${\rm int}(\overline{D})=D$,
$\overline{D}$ is compact (all in the strong topology). These automatically exclude the set $E$ that was constructed in definition 6.5. However, we 
can modify this construction to get at least an open set.

\begin{proposition}
Let $\Gamma$ be any family of holomorphic local homeomorphisms $F\,:\,\mathbb{C}^2\rightarrow\mathbb{C}^2$. Let $U$ be any open subset of $\mathbb{}^2$ with a
smooth boundary that contains the compact $E$. Then the open set $U-E$ is a characteristic set of $\Gamma$.
\end{proposition}
\noindent
{\bf Proof.} \\
Since $E$ can not be mapped in the smooth $\partial U$ by an holomorphic local homeomorphism, we have for any $F_1,F_2\in\Gamma$ for which
$F_1(U-E)=F_2(U-E)$ that also $F_1(E)=F_2(E)$. Now the result follows by Proposition 6.6. $\qed $ \\

\begin{remark}
We note that if $\overline{U}$ is a compact then $U-E$ satisfies, at least the requirement $\overline{U-E}$ is compact. However, the "no slit" 
condition ${\rm int}(\overline{U-E})={\rm int}(\overline{U})\ne U-E$ fails.
\end{remark}
\noindent
Now that we gained some experience with the topological construction of $E$ we are going to make one more step and fix its
shortcomings that were mentioned above. We need to construct a domain $D$ of $\mathbb{C}^2$ which has the following three properties: \\
1) ${\rm int}(\overline{D})=D$ relative to the complex topology. \\
2) $\overline{D}$ is a compact subset of $\mathbb{C}^2$ relative to the strong topology. \\
3) $\forall\,G_1,G_2\in {\rm et}(\mathbb{C}^2)$, $G_1(D)=G_2(D)\Leftrightarrow G_1\equiv G_2$. \\
\\
(The complex topology and the strong topology are the same). Our construction will be a modification of the construction
of the domain that was constructed in Proposition 6.8. We start by modifying the notion of an $m$-star that was introduced
in Definition 6.4.

\begin{definition}
Let $m$ be a natural number and $\alpha\in\mathbb{C}^n$. A thick $m$-star at $\alpha$ is a union of $2m$ triangles, so
that any pair intersect exactly at one vertex, and this vertex (that is common to all the $2m$ triangles) is $\alpha$.
\end{definition}

\begin{definition}
Let $E$ the construction of Definition 6.5 that uses thick $m$-stars.
\end{definition}

\begin{proposition}
Let $\Gamma$ be any family of holomorphic local homeomorphisms $F\,:\,\mathbb{C}^2\rightarrow\mathbb{C}^2$. Then $E$ is
a characteristic set of $\Gamma$.
\end{proposition}
\noindent
{\bf Proof.} \\
The proof is the same word-by-word as that of Proposition 6.6 where we replace $k$-star $S_k$ by thick $k$-star $S_k$. $\qed $ \\
\\
We finally obtain our construction.
\begin{proposition}
Let $\Gamma$ be any family of holomorphic local homeomorphisms $F\,:\,\mathbb{C}^2\rightarrow\mathbb{C}^2$. Let $B(0,R)$ be
an open ball centered at $0$ with a radius $R$ large enough so that $E\subset B(0,R)$ (where $E$ is the set in Proposition 6.12).
Then the domain $D=B(0,R)-E$ is a characteristic set of $\Gamma$.
\end{proposition}
\noindent
{\bf Proof.} \\
The proof is the same as that of Proposition 6.8 where we replace $k$-star $S_k$ by thick $k$-star $S_k$. $\qed $ \\

\section{Injectivity of the left composition operator $L_F$}
We would like our natural mappings: the right mapping $R_F$, and the left mapping $L_F$ to be say
bi-Lipschitz with respect to the metric $\rho_D$ (that reflects the fact that our mappings, ${\rm et}(\mathbb{C}^2)$ satisfy
the Jacobian Condition). Considering first the right mapping $R_F$, it would mean that given three \'etale mappings
$G_1,G_2,F\in {\rm et}(\mathbb{C}^2)$ and a characteristic set $D$ of ${\rm et}(\mathbb{C}^2)$ we need to compare the volume of
$G_1(D)\Delta G_2(D)$ (multiplicities of $G_1$ and of $G_2$ are taken into account) with the volume of the $R_F$ deformed
set, $(G_1\circ F)(D)\Delta (G_2\circ F)(D)$. A short reflection shows that the two volumes are not comparable (in the sense
of bi-Lipschitz). The situation is completely different when we replace the right mapping, $R_F$ by the left mapping, $L_F$.
For example we have the following,

\begin{proposition}
$\forall\,F\in {\rm Aut}(\mathbb{C}^2)$ the mapping $L_F$ is an isometry of the metric space $({\rm et}(\mathbb{C}^2),\rho_D)$.
\end{proposition}
\noindent
{\bf Proof.} \\
For any two mappings $G_1$ and $G_2$ in ${\rm et}(\mathbb{C}^2)$ we need to compare $\rho_D(G_1,G_2)$ with $\rho_D(F\circ G_1,F\circ G_2)$.
We have (using our assumption on $F$), 
$$
(F\circ G_1)(D)\Delta (F\circ G_2)(D)=F\left(G_1(D)\Delta G_2(D)\right).
$$
Since $F$ is also (globally) volume preserving we have,
$$
{\rm the\,volume\,of}\,F\left(G_1(D)\Delta G_2(D)\right)={\rm the\,volume\,of}\,\left(G_1(D)\Delta G_2(D)\right).
$$
This proves that $\rho_D(G_1,G_2)=\rho_D(F\circ G_1,F\circ G_2)$. $\qed $ \\
\\
We now drop the restrictive assumption that $F\in {\rm Aut}(\mathbb{C}^2)$. Thus we merely have $F\in {\rm et}(\mathbb{C}^2)$ and we
still want to compare $\rho_D(G_1,G_2)$ with $\rho_D(F\circ G_1,F\circ G_2)$, for any pair $G_1,G_2\in {\rm et}(\mathbb{C}^2)$. We only know
that $F$ is a local diffeomorphism of $\mathbb{C}^2$ and (by the Jacobian Condition) that it preserves (locally) the volume. In this case
the geometrical degree of $F$, $d_F$ can be larger than $1$. We have the identity $d_F=|F^{-1}(\{(a,b)\})|$ which holds generically (in the
Zariski sense) in $(a,b)\in\mathbb{C}^2$. Hence the (complex) dimension of the set $\{(a,b)\in\mathbb{C}^2\,|\,|F^{-1}(a,b)|<d_F\}$ is at
most $1$. The Jacobian Condition $\det J_F\equiv 1$ implies (as we noticed before) that $F$ preserves volume taking into account the
multiplicity. The multiplicity is a result of the possibility that $F$ is not injective and hence the deformation of the characteristic 
set $D$ by $F$ convolves (i.e. might overlap at certain locations). However, this overlapping is bounded above by $d_F$. So if $A\subseteq\mathbb{C}^2$
is a measurable subset of $\mathbb{C}^2$ and we compare the volume of $A$ with the volume of its image $F(A)$, then,
$$
{\rm the\,volume\,of}\,F(A)\le\,{\rm the\,volume\,of}\,A\le\,d_F\cdot\{{\rm the\,volume\,of}\,F(A)\}.
$$
This can be rewritten as follows,
$$
\frac{1}{d_F}\cdot\{{\rm the\,volume\,of}\,A\}\le\,{\rm the\,volume\,of}\,F(A)\le\,{\rm the\,volume\,of}\,A.
$$
This is the place to emphasize also the following conclusion (that follows by the generic identity $d_F=|F^{-1}(\{(a,b)\})|$), namely
$$
\lim_{A\rightarrow\mathbb{C}^2}\frac{{\rm the\,volume\,of}\,F(A)}{{\rm the\,volume\,of}\,A}=\frac{1}{d_F},
$$
provided that the set $A$ tends to cover the whole of the complex space $\mathbb{C}^2$ in an appropriate manner. To better understand
why the quotient tends to the lower limit $1/d_F$ rather than to any number in the interval $[1/d_F,1]$ (if at all) we recall that
our mapping belongs to ${\rm et}(\mathbb{C}^2)$ and so is a polynomial \'etale mapping. So any point $(a,b)\in\mathbb{C}^2$ for which
$|F^{-1}(a,b)|<d_F$ is an asymptotic value of $F$ and hence belongs to the curve $A_F$ which is the asymptotic variety of $F$. In other words
the identity $d_F=|F^{-1}(a,b)|$ is satisfied exactly on the semi algebraic set $\mathbb{C}^2-A_F$ which is the complement of an algebraic
curve. We now state and prove the main result of this paper,

\begin{theorem}
Let $F,G_1,G_2\in {\rm et}(\mathbb{C}^2)$. Then we have: \\
(i) $\rho_D(F\circ G_1,F\circ G_2)\le \rho_D(G_1,G_2)$. \\
(ii) Suppose that $D$ is a family of characteristic sets of ${\rm et}(\mathbb{C}^2)$ such that $D\rightarrow\mathbb{C}^2$, then
$\forall\,\epsilon>0$ we have,
$$
\left(\frac{1}{d_F}-\epsilon\right)\cdot\rho_D(G_1,G_2)\le\rho_D(F\circ G_1,F\circ G_2)
$$
for $D$ large enough. \\
(iii)  Under the assumptions in (ii) we have:
$$
\lim_{D\rightarrow\mathbb{C}^2}\frac{\rho_D(F\circ G_1,F\circ G_2)}{\rho_D(G_1,G_2)}=\frac{1}{d_F}.
$$
In particular, the left mapping $L_F: {\rm et}(\mathbb{C}^2)\rightarrow {\rm et}(\mathbb{C}^2)$, $L_F(G)=F\circ G$,
is a bi-Lipschitz self-mapping of the metric space $({\rm et}(\mathbb{C}^2),\rho_D)$ with the constants $1/d_F\le 1$.
\end{theorem}
\noindent
{\bf Proof.} \\
(i) $x\in (F\circ G_1)(D)\Delta (F\circ G_2)(D)\Rightarrow \exists\,y\in G_j(D),\,j=1\,{\rm or}\,2$ such that $x=F(y)$ and
$x\not\in (F\circ G_{3-j})(D)$. By $x\not\in (F\circ G_{3-j})(D)$ it follows that $y\not\in G_{3-j}(D)$ and so $y\in G_1(D)\Delta G_2(D)$
and $x=F(y)\in F(G_1(D)\Delta G_2(D))$. Hence $(F\circ G_1)(D)\Delta (F\circ G_2)(D)\subseteq F(G_1(D)\Delta G_2(D))$, so
${\rm vol}((F\circ G_1)(D)\Delta (F\circ G_2)(D))\le{\rm vol}(F(G_1(D)\Delta G_2(D))$, and finally $\rho_D(F\circ G_1,F\circ G_2)\le
\rho_D(G_1,G_2)$. \\
(ii) and (iii). Here the proof is not just set theoretic. We will elaborate more in the remark that follows this proof.
We recall that $F,G_1,G_2\in {\rm et}(\mathbb{C}^2)$. This implies that $\forall\,(\alpha,\beta)\in\mathbb{C}^2$ we have
$|F^{-1}(\alpha,\beta)|\le [\mathbb{C}(X,Y):\mathbb{C}(F)]$, the extension degree of $F$ see \cite{e}. This is the so called Fiber Theorem for
\'etale mappings. Moreover the image is co-finite, $|\mathbb{C}^2-F(\mathbb{C}^2)|<\infty$, \cite{e}.  Also $F$ has a finite
set of exactly $d_F$ maximal domains $\{\Omega_1,\ldots,\Omega_{d_F}\}$. This means that $F$ is injective on each maximal domain
$\Omega_j$, and $i\ne j\Rightarrow \Omega_i\cap\Omega_j=\emptyset$, and $\mathbb{C}^2=\bigcup_{j=1}^{d_F}\overline{\Omega_j}$ and the
boundaries $\partial\Omega_j$ are piecewise smooth (even piecewise analytic). For the theory of maximal domains of entire functions in
one complex variable see \cite{rp1}, and for that theory for meromorphic functions in one complex variable see \cite{sh1,sh2}. Here
we use only basic facts of the theory which are valid also for more than complex variable. If $D$ is a family of characteristic sets of
${\rm et}(\mathbb{C}^2)$ such that $D\rightarrow\mathbb{C}^2$, then by the above $G_1(D),\,G_2(D)\rightarrow\mathbb{C}^2-A$, where $A$
is a finite set, and if $G_1\not\equiv G_2$ then we have the identity,
$$
F(G_1(D)\Delta G_2(D))-(F\circ G_1)(D)\Delta (F\circ G_2(D))=
$$
$$
=\{x=F(y)=F(z)|\,y\in G_1(D)-G_2(D)\,\wedge\,z\in G_2(D)-G_1(D)\}.
$$
Recalling that $(F\circ G_1)(D)\Delta (F\circ G_2(D))\subseteq F(G_1(D)\Delta G_2(D))$) we write the last identity as follows,
$$
F(G_1(D)\Delta G_2(D))=(F\circ G_1)(D)\Delta (F\circ G_2(D))\cup
$$
$$
\cup\{x=F(y)=F(z)|\,y\in G_1(D)-G_2(D)\,\wedge\,z\in G_2(D)-G_1(D)\}.
$$
Taking any two points $y\in G_1(D)-G_2(D)$ and $z\in G_2(D)-G_1(D)$ (as in the defining equation of the set on the right hand side in the last identity), 
we note that there are $i\ne j$, $1\le i,j\le d_F$ such that $y\in\Omega_i\,\wedge\,z\in\Omega_j$ (for $F(y)=F(z)$!).
For $\widetilde{D}$ a large enough characteristic set of ${\rm et}(\mathbb{C}^2)$, we will have $z\in G_1(\widetilde{D})$ and
$y\in G_2(\widetilde{D})$ and so $y,z\in G_1(\widetilde{D})\cap G_2(\widetilde{D})$ (since $G_1(D),G_2(D)\rightarrow \mathbb{C}^2-\{{\rm a\,finite\,set}\}$).
Hence $F(G_1(\widetilde{D})\Delta G_2(\widetilde{D}))-(F\circ G_1)(\widetilde{D})\Delta (F\circ G_2)(\widetilde{D})$ will not include the point $x$.
We conclude that if $y$ and $z$ are $F$-equivalent ($F(y)=F(z)$) then $x=F(y)=F(z)$ will not belong to $F(G_1(D)\Delta G_2(D))-(F\circ G_1)(D)
\Delta (F\circ G_2)(D)$ for large enough $D$. We obtain the following crude estimate:
$$
{\rm vol}(\{x=F(y)=F(z)|\,y\in G_1(D)-G_2(D)\,\wedge\,z\in G_2(D)-G_1(D)\})=
$$
$$
=o({\rm vol}((F\circ G_1)(D)\Delta (F\circ G_2(D)))).
$$
One can think of $D$ as a large open ball centered at the origin of $\mathbb{R}^4$, $D\approx B(R)$ and with the radius $R$ and look at the images of the two
polynomial \'etale mappings $F\circ G_1)(B(R))$ and $(F\circ G_2)(B(R))$ and compare the volume of $(F\circ G_1)(B(R))\Delta (F\circ G_2)(B(R))$ which is of the
order of magnitude $R^{4d}$, where $d$ depends on the algebraic degrees of $F\circ G_1$ and $F\circ G_2$, with the volume of the set
in the left hand side of the last equation. Similar estimates are used in the theory of covering surfaces by Ahlfors, see \cite{h}, chapter 5.We conclude that,
$$
\lim_{D\rightarrow\mathbb{C}^2}\frac{{\rm vol}(F(G_1(D)\Delta G_2(D))}{{\rm vol}((F\circ G_1)(D)\Delta(F\circ G_2)(D)))}=1.
$$
Hence
$$
\lim_{D\rightarrow\mathbb{C}^2}\frac{\rho_D(F\circ G_1,F\circ G_2)}{\rho_D(G_1,G_2)}=\lim_{D\rightarrow\mathbb{C}^2}
\frac{{\rm vol}((F\circ G_1)(D)\Delta(F\circ G_2)(D))}{{\rm vol}(G_1(D)\Delta G_2(D))}=
$$
$$
=\lim_{D\rightarrow\mathbb{C}^2}\frac{{\rm vol}((F\circ G_1)(D)\Delta(F\circ G_2(D))}{{\rm vol}(F(G_1(D)\Delta G_2(D))}\cdot
\frac{{\rm vol}(F(G_1(D)\Delta G_2(D))}{{\rm vol}(G_1(D)\Delta G_2(D))}=
$$
$$
=1\cdot\frac{1}{d_F}=\frac{1}{d_F}.\,\,\,\,\,\qed
$$

\begin{remark}
The facts we used in proving (ii) and (iii) for \'etale mappings are in fact true in any dimension $n$, i.e. in $\mathbb{C}^n$.
In dimension $n=2$ it turns out that the co-dimension of the image of the mapping is $0$ and in fact the co-image is
a finite set. Also the fibers are finite and have a uniform bound on their cardinality (one can get a less tight uniform
bound by the Bezout Theorem). Here are few well known facts (which one can find in Hartshorne's book on Algebraic Geometry, \cite{hart}). \\
1) The following two conditions are equivalent: \\
a. The Jacobian Condition: the determinant $\det J_F$ is a non-zero constant. \\
b. The map $F^{*}$ is \'etale (in standard sense of algebraic geometry). In particular it is flat. \\
Let $F^{*}:\,Y\rightarrow X$ be \'etale. Let $X^{im}:=F^{*}(Y)\subseteq X$. \\
2) For every prime ideal $\wp\subseteq A$ ($X={\rm spec}(A)$), with residue field $k(\wp)$ the ring $B\otimes_A k(\wp)$is finite over
$k(\wp)$ ($Y={\rm spec}(B)$). \\
3) $F^{*}$ is a quasi-finite mapping. \\
4) The set $X^{im}$ is open in $X$. \\
5) For every point $x\in X(\mathbb{C})$ the fiber $(F^{*})^{-1}(x)$ is a finite subset of $Y(\mathbb{C})$. \\
6)The ring homomorphism $A\rightarrow B$ is injective, and the induced field extension $K\rightarrow L$ is finite. \\
7) There is a non-empty open subset $X^{fin}\subseteq X^{im}$ such that on letting $Y^{fin}:=(F^{*})^{-1}(X^{fin})\subseteq Y$,
the map of schemes $F^{*}|_{Y^{fin}}:\,Y^{fin}\rightarrow X^{fin}$ is finite. For any point $x\in X^{fin}(\mathbb{C})$ we have the equality
$d_x=d_{F^{*}}$ the geometrical degree of $F^{*}$. \\
8) The dimension of the set $Z:=X-X^{im}$ is at most $n-2$. \\
9) If $X^{im}=X-Z$ is affine, then $Z=\emptyset$ and $X^{im}=X$. \\
Let $X_{cl}$ be the topological space which is the set $X(\mathbb{C})\cong\mathbb{C}^n$ given the classical topology. Similarly for $Y_{cl}$.
The map of schemes $F^{*}:\,Y\rightarrow X$ induces a map of topological spaces $F_{cl}:\,Y_{cl}\rightarrow X_{cl}$ ($F_{cl}=f^{*}|_{Y(\mathbb{C})}$). \\
10) The map $F_{cl}:\,Y_{cl}\rightarrow X_{cl}$ is a local homeomorphism. \\
\end{remark}
\noindent
An immediate conclusion from Theorem 7.2 is the following,

\begin{corollary}
$\forall\,F\in {\rm et}(\mathbb{C}^2)$ the left mapping $L_F\,:\,{\rm et}(\mathbb{C}^2)\rightarrow L_F({\rm et}(\mathbb{C}^2))$, 
$L_F(G)=F\circ G$ is an injective mapping.
\end{corollary}

\section{Extending the notion of geometrical degree}
In this section we will outline the fact that some of the notions and results that are related to geometrical degree of an \'etale
mapping originate, in fact, in the more basic topological spaces (no algebraic or holomorphic structure is needed). We will skip most of the proofs
(that are elementary).

\begin{definition}
Let $X$ be a topological space. The semigroup of all the continuous mappings, $F:\,X\rightarrow X$, will be denoted by $C(X)$. Here, as usual,
the binary operation is composition of mappings.
\end{definition}

\begin{proposition}
{\rm\bf (1)} Let $X$ be a topological space, $F\in C(X)$ and $\overline{F(X)}=X$. Then $R_F:\,C(X)\rightarrow C(X)$ is injective. \\
{\rm\bf (2)} Let $X$ be a topological space, $F\in C(X)$ and $R_F:\,C(X)\rightarrow C(X)$ is injective. Then for any 
$G,H\in C(X)$, the property $G|_{F(X)}=H|_{F(X)}$ implies that $G\equiv H$, i.e. any $G\in C(X)$ is determined by its restriction $G|_{F(X)}$. \\
{\rm\bf (3)} Let $X$ be a topological space that has the following property: For any closed $C\subseteq X$ and any point $x\in X-C$ there
exist two continuous mappings $F,G:\,X\rightarrow X$ such that $F|_{C}=G|_{C}$ but $F(x)\ne G(x)$. \\
Let $F\in C(X)$ be such that $R_F:\,C(X)\rightarrow C(X)$ is injective. Then $\overline{F(X)}=X$.
\end{proposition}

\begin{remark} Proposition 3.2 follows. \end{remark}
\noindent
We are ready to discuss the notion of the geometrical degree, $d_F$, of appropriate mappings in $C(X)$.

\begin{lemma}
Let $X$ be a topological space, $F\in C(X)$ and $n\in\mathbb{Z}^+\cup\{ 0\}$. If $F$ is open then the set
$B_n=\{x\in X\,|\,|F^{-1}(x)|\le n\}$ is closed.
\end{lemma}

\begin{corollary}
Let $X$ be a topological space, $F\in C(X)$, $F$ open and $n\in\mathbb{Z}^+$. Then we have: \\
{\rm\bf (1)} $\{x\in X\,|\,|F^{-1}(x)|\ne n\}=B_{n-1}\cup (X-B_n)$, the union of a closed set and an open set. \\
{\rm\bf (2)} $\{x\in X\,|\,|F^{-1}(x)|=n\}=B_n-B_{n-1}=B_n\cap(X-B_{n-1})$, the intersection of a closed set 
and an open set. \\
{\rm\bf (3)} If $d_F=\max\{|F^{-1}(x)|\,|\,x\in X\}<\infty$ exists, then $\{x\in X\,|\,|F^{-1}(x)|\ne d_F\}$ is 
a closed set.
\end{corollary}

\begin{remark} $B_n\subseteq B_{n+1}$. \end{remark}

\begin{definition}
Let $X$ be a topological space, $F\in C(X)$, $F$ open and the maximum $d_F=\max\{ |F^{-1}(x)|\,|\,x\in X\}<\infty$ exists.
Then we call $d_F$ the geometrical degree of $F$.
\end{definition}

\begin{example}
If $X=\mathbb{C}^2$ with the complex topology and $F\in {\rm et}(\mathbb{C}^2)$ then we  know that $d_F$ exists. We also know
that the set $A_F=\{ x\in\mathbb{C}^2\,|\,|F^{-1}(x)|<d_F\}$ is a plane algebraic curve (possibly empty). Thus it is
closed in $\mathbb{C}^2$. Moreover it is also small because $\dim A_F<2=\dim\mathbb{C}^2$.
\end{example}
\noindent
We need one more property to hold for our mappings, namely, that the fiber size will generically be $d_F$, i.e., that the set
of all $x\in X$ for which $d_F=|F^{-1}(x)|$ will be a large set measured in the topology of $X$. This leads us to:

\begin{definition}
Let $X$ be a topological space. We will denote by $E(X)$ the set of all the mappings $F:\,X\rightarrow X$ that have
the following properties: \\
{\rm\bf (1)} $F\in C(X)$. \\
{\rm\bf (2)} $F$ is open. \\
{\rm\bf (3)} The maximum $d_F=\max\{ |F^{-1}(x)|\,|\,x\in X\}<\infty$ exists. \\
{\rm\bf (4)} $\overline{X-B_{d_F-1}}=X$.
\end{definition}

\begin{proposition}
Let $X$ be an Hausdorff space. Then: \\
{\rm\bf (1)} $E(X)$ is a semigroup with an identity (where the binary operation is composition of mappings). 
In fact ${\rm Aut}(X)\subseteq E(X)$. \\
{\rm\bf (2)} $\forall\,F,G\in E(X)$, $d_{F\circ G}=d_F\cdot d_G$.
\end{proposition}
\noindent
{\bf Proof.} \\
Checking that $E(X)$ is closed for composition: $F,G\in C(X)\Rightarrow F\circ G\in C(X)$. Also $F,G$ open $\Rightarrow$ $F\circ G$
is open. If $x\in X$ then $|F^{-1}(x)|\le d_F$ and $\forall\,y\in F^{-1}(x)$ we have $|G^{-1}(y)|\le d_G$. Since
$(F\circ G)^{-1}(x)=G^{-1}(F^{-1}(x))$ it follows that $d_{F\circ G}$ exists and, in fact, that $d_{F\circ G}\le d_F\cdot d_G$. This
gives the first three properties in the last definition. We need to check that $\overline{X-B_{d_{F\circ G}-1}}=X$. For that matter
it will be convenient to denote $S_F=\{x\in X\,|\,|F^{-1}(x)<d_F\}$. Let $x\in X-S_G=X-B_{d_G-1}$. This set is open and 
dense in $X$. By our definitions $G^{-1}(x)=\{a_1,a_2,\ldots,a_{d_G}\}$ a finite set of exactly $d_G$ points. Since $X$ is Hausdorff
we can find $d_G$ open neighborhoods $V_{a_1},V_{a_2},\ldots,V_{a_{d_G}}$ of $a_1,a_2,\ldots,a_{d_G}$ respectively which are 
pairwise disjoint. The images $G(V_{a_1}),G(V_{a_2}),\ldots,G(V_{a_{d_G}})$ are open (since $G$ is open) neighborhoods of the
point $x=G(a_1)=G(a_2)=\ldots=G(a_{d_g})$. Let us take the intersection
$$
V=\bigcap_{j=1}^{d_G}G(V_{a_j}).
$$
Then $V$ is an open neighborhood of $x$ and $G^{-1}(V)$ consists of $d_G$ neighborhoods $V_1,V_2,\ldots,V_{d_G}$ of $a_1,a_2,\ldots,a_{d_G}$.
We define $U_j=V_j-S_G$, $j=1,2,\ldots,d_G$. Then since $S_G$ is closed and $S_G^c$ is open and dense in $X$, the $U_j$'s are open and dense subsets
of the $V_j$'s. Since $G$ is continuous and $G(V_j)=V$ it follows that $G(U_j)$ are open and dense in $V$. We note that
each point $y\in U_j$ is such that $|F^{-1}(y)|=d_F$, because $S_F$ is disjoint of $U_j$. Thus the point $x\in X-S_G$ has a neighborhood 
$V$ and $d_G$ open and dense subsets $G(U_j)$ of $V$ such that each point 
$$
x'\in\bigcap_{j=1}^{d_G} G(U_j),
$$
(this set is still open and dense in $V$) is such that $|G^{-1}(x')|=d_G$ and each $a\in G^{-1}(x')$ is such that $|F^{-1}(a)|=d_F$. Hence
$|(F\circ G)^{-1}(x')|=d_F\cdot d_G$. This proves both that $\overline{X-B_{d_{F\circ G}-1}}=X$ and that $d_{F\circ G}=d_F\cdot d_G$. $\qed $ \\

\begin{corollary}
Let $X$ be Hausdorff and $\mathcal{F}$ a semigroup of mappings $F:\,X\rightarrow X$ which are continuous and open and suppose
that there is an absolute constant $c\in\mathbb{Z}^+$ such that $\forall\,F\in\mathcal{F}$, $d_F\le c$. Then the sets $X-B_{d_F-1}$, $F\in\mathcal{F}$
can not be dense in $X$, unless $\forall\,F\in\mathcal{F}$, $d_F=1$.
\end{corollary}

\begin{corollary}
Let $X$ be Hausdorff and $\mathcal{F}$ a family of mappings $F:\,X\rightarrow X$ which are continuous and open and satisfy $\overline{X-B_{d_{F\circ G}-1}}=X$,
$\forall\,F\in\mathcal{F}$ (in particular $\forall\,F\in\mathcal{F}$, $d_F<\infty$). Let $S=<\mathcal{F}>$ be the semigroup generated by $\mathcal{F}$
(composition of mappings is the binary operation). Then $\forall\,F,G\in S$, $d_{F\circ G}=d_F\cdot d_G$.
\end{corollary}

\begin{definition}
Let $X$ be an Hausdorff space. A mapping $F\in E(X)$ is called a composite mapping if $\exists\,G,H\in E(X)-{\rm Aut}(X)$ such that $F=G\circ H$.
A mapping $A\in E(X)-{\rm Aut}(X)$ is a prime mapping if it is not composite. This is equivalent to: if $A=B\circ C$ for some $B,C\in E(X)$,
then $B\in {\rm Aut}(X)\,\vee\,C\in {\rm Aut}(X)$. The subset of $E(X)$ of all the prime mappings will be denoted by $E_p(X)$.
Thus the set of all the composite mappings is $E(X)-E_p(X)$.
\end{definition}

\begin{proposition}
$\forall\,F\in E(X)-E_p(X)$, $d_F$ is not a prime integer. Equivalently, $\forall\,F\in E(X)$, $d_F$ is a prime integer $\Rightarrow\,F\in E_p(X)$.
\end{proposition}

\begin{theorem}
{\rm\bf (1)} If $E(X)-{\rm Aut}(X)\ne\emptyset$, then $E_p(X)\ne\emptyset$. \\
{\rm\bf (2)} $\forall\,F\in E(X)$, $\exists\,k\in\mathbb{Z}^+\cup\{ 0\}$, $\exists\,A_0\in {\rm Aut}(X)$, $\exists\,P_1,\ldots,P_k\in E_p(X)$
such that $F=A_0\circ P_1\circ\ldots\circ P_k$.
\end{theorem}
\noindent
Maybe few elementary examples are in place.

\begin{example}
In Definition 8.9 we take $X=\mathbb{C}$ with the complex topology, and (we an abuse of notation) take $\mathbb{C}[T]-\mathbb{C}$ for $E(X)$. Then
Proposition 8.10(2), $d_{F\circ G}=d_F\cdot d_G$, is the elementary fact from algebra that $\forall\,P,Q\in\mathbb{C}[T]$ we have 
$\deg P\circ Q=\deg P\cdot\deg Q$.
\end{example}
\noindent
A second example is given in section 4 of this paper.

\begin{example}
In definition 8.9 we take $X=\mathbb{C}^2$ with the complex topology, and $E(X)={\rm et}(\mathbb{C}^2)$. Then the theory that was
outlined in Proposition 4.2, Proposition 4.4 and Theorem 4.5 is a special case of the above more general topological theory.
\end{example}

\noindent
{\it Ronen Peretz \\
Department of Mathematics \\ Ben Gurion University of the Negev \\
Beer-Sheva , 84105 \\ Israel \\ E-mail: ronenp@math.bgu.ac.il} \\ 
 
\end{document}